# Urania Propitia

Isobel Falconer
July 2021

Eagle-eyed readers of the July 2021 LMS Newsletter will have spotted that Dr A.E.L. Davis has bequeathed a rare copy of *Urania Propitia* by Maria Cunitz (1604/10 – 1664) to the Society. This significant book was published in 1650, and fewer than twenty-five copies are known to exist.

*Urania Propitia* (Beneficent Urania – Urania is the muse of astronomy) is one of the earliest known books by a female mathematician, publishing cutting edge research in her own name. Its long explanatory subtitle immediately focuses us on some of the mathematical issues of the day. It begins, "wonderfully easy astronomical tables, comprehending the essence of the physical hypotheses brought forth by Kepler, satisfying the phenomena by a very easy, brief way of calculating without any mention of logarithms." We see at once that Kepler's laws are new and require explanation and testing against observation, that logarithms are similarly new and viewed as a barrier rather than an aid to calculation, and that there is a great need for "easy" means of astronomical calculation.

Cunitz perceived this last need as arising as much from the multitude of calendar and almanack makers, and astrologers scattered across Europe, as from the elite astronomers that we hear of more often. She made her twin audiences clear: the book contains parallel texts in Latin and German and she explained that "the German Nation abounds in those of abilities suited to astronomical practice although often lacking knowledge of Latin" while the work was also "intended for the universal improvement of the republic of letters," necessitating Latin [1]. She, herself, came to astronomy and mathematics through the practical route. She was born in 1604 or 1610 in Schweidnitz, Silesia (now Świdnica, Poland) [1]. Both her parents were well educated. Her father, Heinrich Cunitz, was a physician, while her mother, Maria Schultz, was the daughter of Anton Schultz, a practical mathematician. Astrology was still an integral aspect of medicine and Heinrich Cunitz was skilled in this field, as was Maria's second husband, Elias von Löwen, whom she married in 1630.

Prompted by von Löwen, Cunitz became an early user of Kepler's *Rudolphine Tables,* but rapidly concluded that they needed radical transformation if they were to be of extensive practical use. Although the *Rudolphine Tables* were more accurate than previous tables, their basis in Kepler's laws and a heliocentric universe presented two new problems of calculation: finding the position angle of a planet for a given time using Kepler's equation, and combining the heliocentric positions of earth and planet to find its geocentric position [2]. Kepler used logarithmic methods for both, computing the necessary interpolations, but the resultant tables took three times longer to use than earlier tables [1].

Kepler expressed the area traced out by a planet since its *aphelion* (a term that he introduced and Cunitz adopted) in degrees. Cunitz now undertook the c.30,000 calculations necessary to re-cast them in terms of time since aphelion, and provide the true motion per day and per hour for interpolation between the entries for days, and simple double-entry tables from which predictions of planetary position and geocentric latitude and longitude could be found by straightforward arithmetic [1]. *Urania Propitia* took Cunitz around 15 years to complete, interrupted by re-locations necessitated by the Thirty Years War and the birth of three children.

Even before the book's publication in 1650, Cunitz's fame was spreading. She acquired a reputation in her day akin to that accorded to Mary Somerville 200 years later – intellectually outstanding (she knew Latin, Greek and five modern languages), but "modest", "gentle" and skilled in the feminine accomplishments of needlework, painting and music – although she *was* accused of neglecting her domestic affairs [1, quoting Herbinius 1657]. Hevelius, a far more elite and famous astronomer, mentioned her in 1646, and began corresponding with her and von Löwen two years later, discussing a wide range of observational phenomena, new methods and instruments, particularly the telescope [3]. In 1655 fire destroyed the couple's library and instruments, leaving few surviving scraps of Cunitz's manuscripts. Maria Cunitz died in 1664, three years after her husband.

*Urania Propitia* was privately printed, and the print run and purchase price are unknown. Moreover, the work of the jobbing astronomers who may have been its main market has been little studied, so it is difficult to judge Cunitz's success. However, it *was* purchased by James Gregory, first Regius Professor of Mathematics at the University of St Andrews, to equip what would have been Britain's first civic Observatory in 1673 – he did not buy a copy of the *Rudolphine Tables* [4].

Dr A.E.L. Davis was, herself, a Kepler scholar, as well as a great promoter of studies of female mathematicians. She compiled the Davis database of female mathematics graduates in the British Isles [5], and donated her wide-ranging collection of books by and about women who worked in mathematical areas to the Society as the "Philippa Fawcett Collection"[6]. *Urania Propitia* is the magnificent and carefully-chosen bequest of an expert in her fields.

**References**

[1] Swerdlow, N. M. "Urania Propitia, Tabulae Rudophinae Faciles Redditae a Maria Cunitia Beneficent Urania, the Adaptation of the Rudolphine Tables by Maria Cunitz." *A Master of Science History*, edited by Jed Z. Buchwald, vol. 30, Springer Netherlands, 2011, pp. 81–121. *DOI.org (Crossref)*, doi:10.1007/978-94-007-2627-7_7.
[2] Gingerich, Owen. "Johannes Kepler and the Rudolphine Tables." *Resonance*, vol. 14, no. 12, Dec. 2009, p. 1223. doi:10.1007/s12045-009-0116-3.
[3] Włodarczyk, Jarosław. "'Peripheral' Astronomy in the Correspondence of Johannes Hevelius: A Case Study of Maria Cunitia and Elias von Löwen." *Kwartalnik Historii Nauki i Techniki*, vol. 2019, no. Issue 1, Mar. 2019, pp. 147–55. World, *www.ejournals.eu*, doi:10.4467/0023589XKHNT.19.009.10117.

[4] Pringle, R.V. "An edited transcript of the '1687' Catalogue in St Andrews University Library covering the years 1687-1704" http://www.nesms.org.uk/rvp/saul/SAULCat1687.html; see also M. Pilar Gil, "Collecting science in St Andrews. A history in context", M.Litt Dissertation, University of St Andrews, 2016.
[5] "The Davis historical archive: Mathematical Women in the British Isles, 1878-1940", https://mathshistory.st-andrews.ac.uk/Davis/
[6] "The Philippa Fawcett Collection" https://www.lms.ac.uk/library/special-collections#fawcett

Translation of the title page of Urania Propitia (from [1] p86)

BENEFICENT URANIA, or wonderfully easy astronomical tables,
comprehending the essence of the physical hypotheses brought forth by Kepler, satisfying the phenomena by a very easy, brief way of calculating without any mention of logarithms; the concisely taught use of which for the period of time, present, past, and future—adding as well a very easy correction of the superior planets, Saturn and Jupiter, to a more accurate computation and improved agreement with heaven—

MARIA CUNITZ imparts to devotees of this science in a two-fold idiom, Latin and vernacular [German]: That is, new and long-desired, easy astronomical tables, through which the movement of all planets, with regard to longitude, latitude, and other phenomena in all moments of time, past, present, and future, is supplied in an
especially handy method, presented here for the benefit of the devotees of science of the German Nation.